\documentclass{article}
\usepackage{amsmath}
\usepackage{amssymb}

\newcommand{\proof}[1]{\textit{Proof. }#1}

\newtheorem{theorem}{Theorem}
\newtheorem{definition}{Definition}

\newtheorem{remark}{Remark}
\newtheorem{lemma}{Lemma}

\title{A change of variables theorem for the \\ multidimensional Riemann integral
}
\author{Zolt\'an Moln\'ar\\
Department of Algebra,\\
Budapest University of Technology
and Economics,\\
H-1111 Budapest, Egri J.~u.~1.\\
mozow@freemail.hu,\\
Ilona Nagy\\
Department of Mathematical Analysis,\\
Budapest University of Technology
and Economics,\\
H-1111 Budapest, Egri J.~u.~1.\\
nagyi@math.bme.hu,
and\\ Tivadar Szil\'agyi\\
Department of Applied Analysis and
Computational Mathematics\\
E\"otv\"os Lor\'and University,\\
H-1518 Budapest, Pf.120, Hungary,\\
sztiv@math.elte.hu}

\usepackage[mathscr]{eucal}
\newcommand{\ir}{\mathop{\mathrm{int}}\nolimits}
\newcommand{\ext}{\mathop{\mathrm{ext}}\nolimits}

\newcommand{\diam}{\mathop{\mathrm{diam}}\nolimits}
\newcommand{\osc}{\mathop{\mathrm{osc}}\nolimits}
\newcommand{\R}{\ensuremath{\mathbb R}}
\newcommand{\J}{\mathscr{J}}
\newcommand{\D}{\mathscr{D}}
\newcommand{\OO}{\mathscr{O}}
\newcommand{\ee}{\varepsilon}
\newcommand{\dd}{\delta}

\newcommand{\bol}{$\Rightarrow$}
\newcommand{\ol}{\overline}

\begin{document}
\maketitle\thispagestyle{empty}
\begin{abstract}
The most general change of variables theorem for the Riemann
integral of functions of a single variable has been published in
1961 (see \cite{Kestelman:61}). In this theorem, the substitution
is made by an `indefinite integral', that is, by a function of the
form \(t\mapsto c+\int_a^tg=:G(t)\) where \(g\) is Riemann integrable on
\([a,b]\) and \(c\) is any constant. We prove a multidimensional
generalization of this theorem for the case where $G$ is injective
-- using the fact that the
Riemann primitives are the same as those Lipschitz functions which
are almost everywhere strongly differentiable in \((a,b)\). We
prove a generalization of Sard's lemma for Lipschitz functions of several
variables that are almost everywhere strongly differentiable,
which enables us to keep all our proofs within the framework of
the Riemannian theory which was our aim.
\end{abstract}

\renewcommand{\thefootnote}{}
\footnotetext{\hspace*{-.51cm}AMS 2000 subject classification: Primary:
26B10; Secondary: 26B12, 26B15 \\ %
}

\section{Introduction}\label{section1}

As far as we know, the following theorem appeared first in
(\cite{Kestelman:61}).

\begin{theorem} \label{Kestelman:61}
If \(g:[a,b]\to\R\) is Riemann integrable, \(c\in\R\), \(\forall  t\in[a,b]\)
 \(G(t):=c+\int_a^tg\), and the function \(f\) is Riemann integrable
 on the range of \(G\), then \((f\circ G)\cdot g\) is Riemann integrable, and
\[\int_{G(a)}^{G(b)}f=\int_a^b(f\circ G)\cdot g\,.\]
\end{theorem}

Notice that the first statement of Theorem \ref{Kestelman:61} is
somewhat surprising because the composition \(f\circ G\) need not
be Riemann integrable even if \(G\) is \(C^{\infty}\) (see
\cite[Example 34 in Chapter 8.]{Gelbaum:65}). Some years later,
D.~Preiss and J.~Uher in \cite{PreissUher:70} proved the converse:
boundedness of $f$ and integrability of $(f\circ G)g$ implies
integrability of $f$. The aim of the present paper is to formulate and prove
a multidimensional version of these two theorems for the case
where $G$ is injective on the interior of its domain
--- with a proof that remains within the framework of the Riemann
theory. In textbooks, the usual assumption on $G$ is that it has a continuously
differentiable extension to an open set that covers the closure of the original
-- Jordan measurable -- domain. Observe that this assumption implies Lipschitz
continuity (see Theorem \ref{locLip}). The starting point to the corresponding
generalization is the fact (which
seems not to be well-known) that a function $G:[a,b]\to\R$ is a
Riemann primitive if and only if it is Lipschitz and almost
everywhere strongly differentiable.

In the next section, after introducing some notations and
terminology, we summarize some well-known facts about Riemann
integrability, and give a basic theorem (\ref{alap}) about the
change of variables with easy proof and 'hard-to-check'
assumptions. In section 3, we investigate the notion 'strong
differentiability' and other auxiliary tools, then in section 4
we prove that for
injective functions $G$ that are Lipschitz and almost everywhere
strongly differentiable, and for properly chosen $g$, conditions
a), b)and c) of Theorem\,\ref{alap} are fulfilled (injectivity
will be assumed only on the difference of $D(G)$ and a set of
Lebesgue measure zero).

\section{Terminology and some basic facts about Riemann
integrability}\label{section2}

For any $H\subset\R^m$, the set of Jordan measurable subsets of
$H$ will be denoted by $\mathscr{J}_H$, in the case $H:=\R^m$ the
subscript will be omitted. The volume or Jordan content of a
Jordan measurable set $X\subset\R^m$ will be denoted by $ V(X)$
and the outer Jordan content of a bounded set $Y$ by $V^*(Y)$.

By a Jordan partition of $X\in\J$ we mean a finite collection of
pairwise non-over\-lap\-ping sets in $\J_X$ the union of which is $X$.
The set of all Jordan partitions of $X\in\J$ will be denoted by
$\Pi(X)$. By the norm of a Jordan partition $\Phi\in\Pi(X)$ we
mean the number $|\Phi|:=\max\{\diam(H)\,:\,H\in\Phi\}$. The lower
sum, upper sum and oscillation sum of a bounded function
$f:X\to\R$ corresponding to the partition $\Phi\in\Pi(X)$ is
defined by $s_f(\Phi):=\sum_{H\in\Phi}\inf f\vert_H V(H)$,
$S_f(\Phi):=\sum_{H\in\Phi}\sup f\vert_H V(H)$, resp.
\[\OO_f(\Phi):=S_f(\Phi)-s_f(\Phi)=\sum\osc_f(H)V(H),\quad\mbox{
where}\quad \]
\[\osc_f(H):=\sup f\vert_H-\inf f\vert_H=
\sup\{|f(y)-f(x)|\,:\,x\in H,\ y\in H\}.\]

\noindent The lower and upper Darboux integral of $f$ is
\[\sideset{}{_X}\int\limits_{\ol{\phantom{J}}\ \ \ } f:=\sup
s_f,\qquad\mbox{ resp.}\qquad \sideset{}{_X}\int\limits^{\ \ -}
f:=\inf S_f.\]
The bounded function $f\colon X\to\R$ is integrable
(with integral $\alpha\in\R$) if its lower and upper Darboux
integrals agree (and are equal to $\alpha$).

\vspace{1mm} By a dotted Jordan partition of $X\in\J$ we mean a
finite set of ordered pairs
\[\eta:=\{(H_1,y_1),\ldots,(H_n,y_n)\}\] such that
$D(\eta):=\{H_1,\ldots,H_n\}\in\Pi(X)$, and $y_i\in H_i$ for
$\,i=1,\ldots,n$. The Riemann sum of the function $f\colon X\to
\R$ corresponding to the dotted Jordan partition $\eta$ is
$\sigma_f(\eta):=\sum_{i=1}^nf(y_i)V(H_i)$.

We will make use of the following well-known statements:

\begin{theorem}[Generalized Darboux Theorem]\label{Darb}
For each $X\in\J$ and for each boun\-ded $f\colon X\to\R$,
\[\lim_{|\Phi|\to0}s_f(\Phi)=\sideset{}{_X}\int\limits_{\ol{\phantom{J}}\ \ \ } f,\qquad
\lim_{|\Phi|\to0}S_f(\Phi)=\sideset{}{_X}\int\limits^{\ \ -} f.\]
\end{theorem}

\begin{theorem}[modified Riemann's condition]\label{oszc}
For each $X\in\J$ and for each bounded $f\colon X\to\R$, integrability of
$f$ is equivalent to the condition $\lim_{|\Phi|\to0}\OO(\Phi)=0$.
\end{theorem}

\begin{theorem}\label{Rdef}
For each $X\in\J$, $\alpha\in\R$ and $f\colon X\to\R$, the
following two statements are equivalent: 1. $f$ is integrable with
integral $\alpha$, 2. $\lim_{|D(\eta)|\to0}\sigma_f(\eta)=\alpha$.
\end{theorem}

The definition of the integral based on Riemann sums can be used
in the matrix-valued case, too. In $\R^{m\times n}$, any metric
induced by a norm can be used. In particular, for each $X\in\J$
and integer $m>1$, a matrix-valued function $h\colon
X\to\R^{m\times m}$ is integrable if and only if all the entries
$h_{ik}\colon X\to\R$ $(i,k=1,\ldots,m)$ are integrable. This fact
will be used in order to simplify the formulation of our last
theorem.

\begin{theorem}\label{addi}
For each $X\in\J$, $H\in\J_X$ and integrable $f\colon X\to\R$, the
restriction $f\vert_H$ is integrable, for each $\Phi\in\Pi(X)$,
$\int_Xf=\sum_{H\in\Phi}\int_Hf$.
\end{theorem}

\begin{theorem}\label{absz}
For each $X\in\J$ and integrable $f\colon X\to\R$, the function
$|f|$ is integrable, and the inequality
$\left|\int_Xf\right|\le\int_X|f|$ holds.
\end{theorem}

\begin{definition}
If $X\subset\R^m$, $\D\subset\J_X$, $\cup\D=X$ and
$\Psi\colon\D\to\R$, then by a density function of $\Psi$ we mean a
function $g\colon X\to\R$ for which integrability of $g|_H$ and
$\Psi(H)=\int_Hg$ hold for each $H\in\D$.
\end{definition}

\begin{theorem}\label{alap}
Let $X\in\J$, $G\colon X\to\R^m$ a Lipschitz function, $Y:=G(X)$,
$f\colon Y\to\R$ bounded and $g\colon X\to\R$ integrable. Suppose
that

\noindent a) for each $H\in\J_X$, $G(H)\in\J$,

\noindent b) for each pair of non-overlapping sets $A\in\J_X$,
$B\in\J_X$, $G(A)$ and $G(B)$ are non-overlapping,

\noindent c) $g$ is a density function of the function $\J_X\ni
H\mapsto V(G(H))$,

\noindent d)$f$ or $(f\circ G)g$ is integrable.

Then both $f$ and $(f\circ G)g$ are integrable and $\int_X(f\circ
G)g=\int_Yf$.
\end{theorem}
\proof{ Let $L>0$ be a Lipschitz constant for $G$, $K>0$ such that
for every $x\in Y$, $|f(x)|\le K$, and use the notation
$\psi:=(f\circ G)g$.

First, suppose that $f$ is integrable and let $\ee$ be a positive
number. We will show that for some $\dd>0$, and for each dotted
partition $\eta$ of $X$ with $|D(\eta)|<\dd$, we have
$|\sigma_{\psi}(\eta)-\int_Yf|<\ee$ (see Theorem\,\ref{Rdef}).
Choose positive numbers $\dd_f$ and $\dd_g$ such that
$\OO_f(\Psi)<\ee/2$ whenever $\Psi\in\Pi(Y)$ and $|\Psi|<\dd_f$,
resp. $\OO_g(\Phi)<\ee/2K$ whenever $\Phi\in\Pi(X)$ and
$|\Phi|<\dd_g$ (see Theorem\,\ref{oszc}). Let
$\{(H_k,y_k)\,:\,k=1,\ldots,n\}$ a dotted partition of $X$ such
that the norm of $\Phi:=\{H_1,\ldots,H_n\}$ is less then
$\min\{\dd_g,\dd_f/L\}=:\dd$. Conditions a) and b) imply that
$\Psi:=\{G(H_1),\ldots,G(H_n)\}$ is a Jordan partition of $Y$,
Lipschitz condition and the definition of $L$ imply that the norm
of this latter partition is smaller then $\dd_f$.
\[\left|\sum_{k=1}^nf(G(y_k))g(y_k)V(H_k)-\int_Yf\right|\,\stackrel{\ref{addi}}{=}\,
\left|\sum_{k=1}^n\big(f(G(y_k))g(y_k)V(H_k)-\int_{G(H_k)}f\big)\right|\]
\[=\left|\sum_{k=1}^n\left[f(G(y_k))\left[g(y_k)V(H_k)-V(G(H_k))\right]+
f(G(y_k))V(G(H_k))-\int_{G(H_k)}f\right]\right|\]
\[\stackrel{c)}{\le}\,\sum_{k=1}^n\left|f(G(y_k))\right|\left|g(y_k)V(H_k)-\int_{H_k}g\right|+
\sum_{k=1}^n\left|f(G(y_k))V(G(H_k))-\int_{G(H_k)}f\right|\]
\[=\sum_{k=1}^n\left|f(G(y_k))\right|\left|\int_{H_k}[g(y_k)-g(x)]\,dx\right|+
\sum_{k=1}^n\left|\int_{G(H_k)}[f(G(y_k))-f(y)]\,dy\right|\]
\[\,\stackrel{\ref{absz}}{\le}\,K\sum_{k=1}^n\int_{H_k}|g(y_k)-g(x)|\,dx+
\sum_{k=1}^n\int_{G(H_k)}|f(G(y_k))-f(y)|\,dy\]
\[\le K\sum_{k=1}^n\osc_g(H_k)V(H_k)+\sum_{k=1}^n\osc_f(G(H_k))V(G(H_k))\]
\[=K\OO_g(\Phi)+\OO_f(\Psi)<\dfrac{\ee}{2}+\dfrac{\ee}{2}.\]

Second, suppose that $\psi$ is integrable, we prove that
$\int_X\psi$ is equal to the upper Darboux integral of $f$. The
proof of the fact that $\int_X\psi$ is equal to the lower integral
of $f$ is completely similar, therefore it will be omitted. Let
$\ee$ be a positive number, we show that $\int_X\psi$ is in the
$\ee$-neighborhood of the upper integral of $f$. According to
Theorem\,\ref{Rdef}, one can choose a $\dd_{\psi}>0$ such that
$|\sigma_{\psi}(\eta)-\int_X\psi|<\ee/4$ holds for every dotted
partition $\eta$ of $X$ with $|D(\eta)|<\dd_{\psi}$, according to
Theorem\,\ref{oszc} -- a $\dd_g>0$ such that $\OO_g(\Phi)<\ee/4K$
holds whenever the norm of $\Phi\in\J_X$ is less then $\dd_g$ and
according to Theorem\,\ref{Darb} -- a $\dd_f>0$ such that
$S_f(\Psi)$ lies in the $\ee/4$-neighborhood of the upper integral
of $f$ whenever the norm of $\Psi\in\J_Y$ is less then $\dd_f$.
Fix a Jordan partition $\Phi=\{H_1,\ldots,H_n\}\in\J_X$ with
\[|\Phi|<\min\{\dd_{\psi},\dd_g,\dd_f/L\}=:\dd,\]
and for each $k=1,\ldots,n$ an element $y_k\in H_k$ such that
\[f(G(y_k))>\sup f\vert_{G(H_k)}-\dfrac{\ee}{4(V(Y)+1)}.\]
Denoting the collection of sets $G(H_k)$ by $\Psi$ and the set of
pairs $(H_k,y_k)$ by $\eta$ (k=1,\ldots,n), we have
\[\int_X\psi-\sideset{}{_Y}\int\limits^{\ \ -} f=\left[\int_X\psi-\sigma_{\psi}(\eta)\right]+
\sum_{k=1}^nf(G(y_k))[g(y_k)V(H_k)-V(G(H_k))]\]
\[+\sum_{k=1}^n[f(G(y_k))-\sup f\vert_{G(H_k)}]V(G(H_k))+
\left[S_f(\Psi)-\sideset{}{_Y}\int\limits^{\ \ -} f\right].\]
$\dd\le\dd_{\psi}$, $\dd\le\dd_g$, the choice of the points $y_k$
and $\dd\le\dd_f$ imply that the absolute value of the first,
second, third, respectively the fourth member on the right hand
side is less then $\ee/4$. As for the second member, this is seen
from the following estimate:
\[\left|\sum_{k=1}^nf(G(y_k))[g(y_k)V(H_k)-V(G(H_k))]\right|\stackrel{c)}{=}
\left|\sum_{k=1}^nf(G(y_k))\int_{H_k}[g(y_k)-g(x)]\,dx\right|\]
\[\le\sum_{k=1}^n|f(G(y_k))|\left|\int_{H_k}\osc_g(H_k)\right|\le
K\OO_g(\Phi).\]
}

\section{Auxiliary tools}\label{section3}

\subsection{Strong differentiability}
\begin{definition}\label{strongdif}
Let $m$ and $ n $ be positive integers, $ U\subset{\bf R}^m $, and $u$
an interior point of $U$. The function $f:U\to{\bf R}^n$ is \emph{strongly differentiable}
at the point $u$, if there exists a linear map $ A:{\bf R}^m\to
{\bf R}^n $ such that
$$ \lim_{(x,y)\to(u,u)}\frac{1}{\|x-y\|}\left[ f(x)-f(y)-A(x-y)\right]\,=\,0,$$
where on $\R^m\times\R^m$ one can use any metric induced by a norm, e.g.  \newline
$d((x,y),(z,w)):=\max\{\|x-z\|,\|y-w\|\}$.
\end{definition}

We make some remarks about this notion.

Some authors use the term `strict differentiability' instead of
`strong diffe\-rentiability.'

In the definition the spaces $ \R^m $ and $ \R^n $ could be
replaced by any normed spaces, but in this case (if the first
space is infinite dimensional) one says `continuous linear'
instead of `linear'.

If $f$ is strongly differentiable at $ u $ then it is differentiable there and
$f'(u)=A $ must hold.

Strong differentiability of $ f $ at $ u $ implies the existence of a neighborhood
of $ u $ on which $ f $ is a Lipschitz function.

If one replaces the assumption on continuous differentiability of
$f$ at the point $ u $ by strong differentiability at $u$ in the
local inverse function theorem (see \cite{Graves:27}, then one can
state existence of a neighborhood $U$ of $u$ such that the
restriction $f|_U$ is \emph{injective}, its range is a
neighborhood of $ f(u) $, the local inverse (the inverse of this
restriction) is strongly differentiable at the point $f(u) $, and
the derivative of the local inverse at $f(u)$ is equal to the
inverse of $f'(u)$. (see for example \cite{Leach:61} or
\cite{Leach:63}). Tho most difficult part of the proof is essentially contained in
the proof of the next theorem (the proof of the fact that $f(u)$ is an interior point of the range of the injective restriction).

\begin{theorem}\label{invfunthe}
Let $u$ be an element of the open set $\Omega\subset\R^m$, suppose
that the function $G:\Omega\to\R^m$ is strongly differentiable at
$u$ and $M:=G'(u)$ is regular. Then for each $\ee\in(0,1)$, there
exists a $\dd>0$ such that
$\forall(x,r)\in\Omega\times(0,+\infty)$,
\[ \ol{B}(x,r)\subset\ol{B}(u,\delta)\;
\Longrightarrow\;\ell_x\ol{B}(x,(1-\ee)r)\subset G(\ol{B}(x,r))
\subset\ell_x\ol{B}(x,(1+\ee)r),\] where $\ell_x:\R^m\to\R^m$
denotes the affine function $z\mapsto G(x)+M(z-x)$. In particular,
$G(u)\in\ir R(G).$
\end{theorem}

\proof{ Define the function $ \varrho:\Omega\times\Omega\to\R^m $
as follows: if $x,z\in\Omega$ and $x=z$ then $\varrho(x,z):=0$,
otherwise
\[\varrho(x,z):=\frac{1}{\|z-x\|}[G(z)-G(x)-M(z-x)].\]
The strong differentiability condition implies that for each
$\ee\in(0,1)$, one can find a $\delta>0$ such that
\[\|\varrho(x,z)\|<\frac{\ee}{\|M^{-1}\|},\qquad\mbox{whenever}\quad
x,z\in\ol{B}(u,\delta).\] Fix a pair $(x,r)$ satisfying the
condition $\ol{B}(x,r)\subset\ol{B}(u,\delta)$, and, in order to
prove the first inclusion, fix an element $y=\ell_x(v)$ with
$\|v-x\|\le(1-\ee)r$ as well. One can apply Banach's fixed point
theorem on the metric subspace $X:=\ol{B}(x,r)$ of $\R^m$ to the
function
\[f:X\to\R^m,\quad\quad z\mapsto z-M^{-1}[G(z)-y]=M^{-1}[y-G(z)+Mz]\]
(of course, the fixed point is a point $z\in\ol{B}(x,r)$ with
$G(z)=y$). Indeed, $f$ maps $X$ into $X$, because for each $z\in
X$ we have
\[\|f(z)-x\|=\|M^{-1}[y-G(x)+G(x)-G(z)-M(x-z)]\|\le\]
\[\|M^{-1}[\ell_x(v)-G(x)]\|+\|M^{-1}\|\|G(x)-G(z)-M(x-z)\|\le\]
\[\|v-x\|+\ee\|x-z\|\le(1-\ee)r+\ee r=r,\]
and it is a contraction with Lipschitz constant $\ee$, because for
each pair $(z,w)\in X\times X$ we have
\[\|f(z)-f(w)\|=\|M^{-1}[G(w)-G(z)-M(w-z)]\|\le
\|M^{-1}\|\frac{\ee}{\|M^{-1}\|}\|w-z\|.\] To prove the second
inclusion, fix an element $v\in\ol{B}(x,r)$ and set
\[w:=M^{-1}[\|v-x\|\varrho(x,v)].\] Then
$\|w\|\le\|M^{-1}\|\cdot\|\varrho(x,v)\|\cdot\|v-x\|$, hence
$\|v+w-x\|\le r+\varepsilon r$, and
\[G(v)=G(x)+M(v-x)+\|v-x\|\varrho(x,v)=G(x)+M(v+w-x)=\ell_x(v+w). \]
To show that $G(u)$ is an interior point of the range, apply the
first inclusion with $\ee:=1/2$, $x:=u$, $r:=\dd$:
\[G(\ol{B}(u,\dd))\supset G(u)+MB(0,\dd/2)\supset
G(u)+B(0,\dd/(2\|M^{-1}\|)).\]
}

\begin{remark}\label{edked}
Let $ m $, $ n $, $ U $ and $ u $ be the same as in Definition
\ref{strongdif}. The function $f:U\to\R^n$ is strongly
differentiable at $ u $ if and only if each component of $ f $ has
this property.\end{remark}
\begin{remark} \label{contofstd}
Let $ m $, $ n $, $ U $ and $ u $ be the same as in Definition
\ref{strongdif}. If the function $ f:U\to\R^n $ is strongly
differentiable at $u$, then, without any assumption on the domain
of $f'$, the function $ x\mapsto f'(x)$ is continuous at $ u $.
For a proof, see \cite{Nijenhuis:74}.\end{remark}
\begin{remark}
Let $ m $, $ n $, $ U $ and $ u $ be the same as in Definition
\ref{strongdif}. If the function $ f:U\to\R^n $ is differentiable
in a neighborhood of $ u $ and $ f' $ is continuous at $ u $ then
$ f $ is strongly differentiable at $ u $. As for the proof: apply
the mean value inequality to the function $ z\mapsto f(z)-f'(u)z $
on the line segment $[x,y]$. \end{remark}

Before our last remark we introduce a definition which is a slight
modification of Nijenhuis' definition (see \cite{Nijenhuis:74}).

\begin{definition}\label{epd}
Let $m$ and $j$ be positive integers, $j\le m$, $ U\subset{\bf
R}^m $, and $u$ an interior point of $U$. The function $f:U\to\R$
is \emph{strongly partially differen\-tiable} with respect to the
$j$-th variable at the point $u$, if there exists a real number
$D^s_jf(u)$ such that for each $ \varepsilon>0 $ there is a $
\delta>0$ with the following property: if $ x,y\in B(u,\delta) $,
$ y_j\neq x_j $, but for all $ i\neq j $ $ y_i=x_i $, then
$$\left|\frac{f(x)-f(y)}{x_j-y_j}-D^s_jf(u)\right|<\varepsilon\,.$$
\end{definition}

\begin{remark}\label{edepd}
Let $U$ be a subset of ${\bf R}^m $, and $u$ an interior point of
$U$. The function $f:U\to\R$ is strongly differentiable at $ u $
if and only if $ f $ is strongly partially differentiable at $ u $
with respect to all its variables. (The proof is an easy exercise,
for the case $m=2$ see \cite{Nijenhuis:74}.)
\end{remark}

\subsection{Cubes, set-functions}
By a (closed) cube we mean a Cartesian product of \(m\) number of
closed one dimensional intervals of equal length, by a
cube-partition of a cube \(Q\) we mean a finite set of pairwise
non-overlapping \emph{cubes}, the union of which is \(Q\).
Analogously, by a dotted cube-partition of a cube \(Q\) we mean a
finite set
\[ \{(Q_1,y^1),\ldots,(Q_n,y^n)\}\]
of ordered pairs, where \(\{Q_1,\ldots ,Q_n\}\) forms a cube-partition of \(Q\) and \(y^i\in Q_i\) for \(i=1,\ldots ,n\).
Equivalently, dotted cube-partitions of a cube \(Q\) can be viewed as functions: a function \(\eta:\mathcal{A}\to Q\)
is a dotted cube-partition of the cube \(Q\), if \(\mathcal{A}\) is a cube-partition of \(Q\) and for each cube
\(I\in \mathcal{A}\), \(\eta _I:=\eta (I)\in I\). In the following, (dotted) partition of a cube will mean always a (dotted) cube-partition.

In the space \(\R^m\) we use the norm \(x\mapsto\max|x_i|=:\|x\|\), therefore the closed balls \(\ol{B}(y,r)\)
are cubes and the open balls $ B(y,r) $ are open cubes.

Let \(\delta\) be a positive valued function defined on a cube \(Q\). A dotted partition \(\eta\) of \(Q\)
is said to be \(\delta\)-fine, if for each \(I\in
D(\eta)\), \(I\subset B(\eta _I,\delta(\eta _I))\).
The following statement will be called `Cousin's lemma':

\begin{lemma}[Cousin's lemma]\label{Couslemma}
For each cube \(Q\) and each function \(\delta:Q\to(0,+\infty)\), \(Q\) has a  \(\delta\)-fine dotted partition.
\end{lemma}

For a proof of this assertion, see for example the proof of
\cite[Lemma 7.3.2]{Pfeffer:93}.

\begin{definition}
A real valued function $\Phi$ defined on a set $\mathcal{C}\subset\J$ will be called
\begin{enumerate}
\item \emph{additive}, if for each $ H\in\mathcal{C} $ and each Jordan
partition $\mathcal{A}\subset\mathcal{C}$ of $H$, $$
\Phi(H)=\sum_{J\in\mathcal{A}}\Phi(J), $$ \item \emph{Lipschitz},
if there exists a nonnegative number $ L $ such that for each $
H\in\mathcal{C} $, $ \Phi(H)\le L\cdot V(H) $,
\end{enumerate}
\end{definition}

\begin{definition}
Suppose that for some $r>0$ and $u\in\R^m $, each closed subcube
of $B(u,r)$ belongs to the domain of the real valued set-function
$\Phi$. Then $\Phi$ is called \emph{differentiable}
(resp.~\emph{strongly differentiable}) at  $u$,  if for some real
number $ \Phi'(u) $ and for each $\omega>0 $ there exists a
$\dd>0$ for which $ u\in I\in D(\Phi)$ (resp.~$I\in D(\Phi)$) and
$I\subset B(u,\dd)$ imply
$$ \left|\frac{\Phi(I)}{V(I)} - \Phi'(u)\right|\,<\,\omega. $$
\end{definition}

Of course, if a cube-function $\Phi$ is differentiable at $u$ then
the number $\Phi'(u)$ in the definition of differentiability is
unique. Note that if $m=1$, $f:[a,b]\to\R$ and $\Phi$ is defined
on the set of closed subintervals of $[a,b]$ by
$[\alpha,\beta]\mapsto f(\beta)-f(\alpha)$, then $ \Phi $ is
additive, $ \Phi $ is Lipschitz if and only if $ f $ is Lipschitz,
$ \Phi $ is (strongly) differentiable at $ u\in(a,b)$ if and only
if $ f $ is (strongly) differentiable there.

\begin{remark}\label{indef}
Another important and well-known example: if $X\in\J$ and $g\colon
X\to\R$ is integrable, then the set-function $\J_X\ni
H\mapsto\int_Hg=:\Psi(H)$ is additive, Lipschitz, and strongly
differentiable at the continuity points $u$ of $g$ with
$\Psi'(u)=g(u)$.
\end{remark}

\subsection{Lipschitz functions and Lipschitz set-functions}

\begin{lemma}\label{Lipkiterj}
If $(M,d)$ is a metric space, $\emptyset\neq X\subset M$ and
$G\colon X\to\R^m$ is a Lipschitz function with Lipschitz constant
$L$, then $G$ has a Lipschitz extension $F\colon M\to\R^m$ with
the same Lipschitz constant $L$.
\end{lemma}
\proof{ Because of our choice of the norm in $\R^m$, the lemma
follows from the special case where $m=1$, which can be applied to
the component functions. But this special case is a well-known
theorem, for a proof see for example \cite[6.6.5 and
6.6.6]{Pfeffer:93}. }

\begin{theorem}\label{VGHLip}
Let $X\in\J$ and $G\colon X\to\R^m$ a Lipschitz function, then the
set-function $\J_X\ni H\mapsto V^{*}(G(H))=:\Psi(H)$ is again
Lipschitz.
\end{theorem}
\proof{ Let $F\colon\R^m\to\R^m$ be an extension of $G$ satisfying the
Lipschitz condition with Lipschitz constant $L$ (see
Lemma\,\ref{Lipkiterj}), so for each cube
$I=\ol{B}(u,r)\subset\R^m$, Lipschitz condition yields
$F(I)\subset\ol{B}(F(u),Lr)$, consequently $V^{*}(F(I))\le
L^mV(I)$. Let $H\in\J_X$, $\ee>0$ and $\{I_1,\ldots,I_n\}$ be a
finite set of cubes such that
\[H\subset\cup_{k=1}^nI_k\quad\mbox{and}\quad\sum_{k=1}^nV(I_k)<V(H)+\ee.\]
$V^{*}$ is monotonic and subadditive, therefore
\[V^{*}(G(H))\le V^{*}\left[\cup_{k=1}^nF(I_k)\right]\le
\sum_{k=1}^nV^{*}(F(I_k))\le L^m\sum_{k=1}^nV(I_k)<L^m(V(H)+\ee)\]
and this gives the inequality $\Psi\le L^mV\vert_{\J_X}$.
}

\begin{theorem}\label{locLip}
If $(X_1,d_1)$ is a compact metric space, $(X_2,d_2)$ a metric
space and $f\colon X_1\to X_2$ is locally Lipschitz in each point
of $X_1$ then $f$ is Lipschitz.
\end{theorem}

\proof{
Using on $X_i\times X_i$ -- for example -- the metric
\[((x,y),(u,v))\mapsto\max\{d_i(x,u),d_i(y,v)\}\qquad(i=1,2),\]
$X_1\times X_1$ is compact, therefore (being a closed subset of
this compact space) the diagonal $\Delta:=\{(x,x)\,:\,x\in X_1\}$
is also compact. This fact and the local Lipschitz condition gives
a positive integer $n$, elements $z_1,\ldots,z_n\in X_1$ and
positive numbers $r_1,\ldots,r_n,L_1,\ldots,L_n$ such that for
each $k=1,\ldots,n$, $f\vert_{B(z_k,r_k)}$ is Lipschitz with
Lipschitz constant $L_k$, and
\[\Delta\subset\cup_{k=1}^nB(z_k,r_k)\times B(z_k,r_k)=:\Gamma.\]
$(X_1\times X_1)\setminus\Gamma$ is again compact, the metrics and
$G$ are continuous, thus the restriction to $(X_1\times
X_1)\setminus\Gamma$ of the function
\[(X_1\times
X_1)\setminus\Delta\ni(x,y)\mapsto\frac{d_2(f(x),f(y))}{d_1(x,y)}=:h(x,y)\]
has an upper bound $L_0$. This implies that
$\max\{L_0,L_1,\ldots,L_n\}$ is an upper bound of $h$.
}

\subsection{Some consequences of Cousin's lemma}

The following lemma is known from several proofs of Sard's lemma
(see for example \cite[proof of Theorem~3.14.]{Spivak:65}).

\begin{lemma}\label{lemmaSard}
Suppose that $A$ is a subset of $\R^m $, $G:A\to\R^m$ is strongly
differentiable at the interior point $u $ of $A$ and let $G'(u)$
be singular. Then for each $\varepsilon>0$ there exists a
$\delta>0$ such that $ B(u,\delta)\subset A$, and the inequality $
V^*(G(I))\le\varepsilon V(I)$ holds for all cubes $I$ covered by
$B(u,\delta)$.
\end{lemma}

Now we prove an interesting version of the so-called Sard's lemma.
Observe that as A.~Sard himself writes in \cite{Sard:42}, the real
valued $ C^1 $-case is due to A.~P.~Morse (see \cite{Morse:39}).If
$G$ is differentiable at an interior point $ x $ of its domain
then the Jacobi matrix of $G$ at $x$ will be denoted by $J_G(x)$.

\begin{theorem} \label{genofSard}
Suppose that $\Omega$ is an open subset of $\R^m$,
$K\subset\Omega$ is a set of Lebesgue measure $0 $ and
$G:\Omega\to\R^m$ is a Lipschitz function which is strongly
differentiable at every points of $\Omega\setminus K$. Then the
image under $ G $ of the set
\[\{x\in\Omega\setminus K\,:\,J_G(x)\mbox{ is singular }\}  \]
is of Lebesgue measure $ 0 $.
\end{theorem}

\proof{ $\Omega$ is a countable union of cubes, so it is enough to
prove that for any cube $Q\subset\Omega $, the image under $ G $
of the set
\[S:=\{x\in Q\setminus K\,:\,J_G(x)\mbox{ is singular }\}  \]
has Jordan content $ 0 $. Set
\[R:=\{x\in Q\setminus K\,:\,J_G(x)\mbox{ is regular }\},\quad
\mbox{and}\quad T:=Q\cap K.\] First, observe that $ R\subset\ext
S$. Indeed, $J_G$ is continuous at every point of $R$ (see Remark
\ref{contofstd}) and so is the function $\det:\R^{m\times
m}\to\R$, thus a neighborhood of a point of $ R $ in which for
every $ x $ we have $\det J_G(x)\neq0 $ cannot intersect $ S $.
Second, observe that if $ L $ is a Lipschitz constant for $ G $
then for any cube $I:=\ol{B}(u,r)\subset\Omega$ we have $
G(I)\subset\ol{B}(G(u),Lr)$, so $ V^*(G(I))\le L^mV(I) $. In order
to apply Cousin's lemma, we define a positive valued function $
\delta $ on $ Q $. Let $ \varepsilon $ be a positive number, fix a
countable set $ \mathcal{I} $ of open intervals with the sum of
volumes being less than $\ee/2L^m$, the union of which covers the
set $ T $. If $ u\in R $ then let $\delta(u)>0$ be such that
$B(u,\delta(u))\cap S=\emptyset$. For each $ u\in S$, using the
previous lemma, select a $ \delta(u)>0 $ such that $
V^*(G(I))\le\ee V(I)/2V(Q) $ holds for every cube $ I $ satisfying
conditions $ u\in I\subset B(u,\delta(u)) $. For each $ u\in T $,
select first a $ J_u\in\mathcal{I} $ that contains $ u $ and then
a $ \delta(u)>0 $ such that $B(u,\delta(u))\subset J_u $. Fix a $
\delta $-fine dotted partition $ \eta $ of $ Q $ and write its
domain in the form $ \mathcal{A}\cup\mathcal{B}\cup\mathcal{C} $
where for a cube $ I\in D(\eta) $, $ I\in\mathcal{A} $ means $
\eta_I\in R $, $ I\in\mathcal{B}$ means $\eta_I\in S$ and
$I\in\mathcal{C}$ means $\eta_I\in T$. As $\eta$ is $ \delta
$-fine, each $ I\in\mathcal{A} $ is disjoint from $ S $, hence
\[G(S)\subset\bigcup_{I\in\mathcal{B}\cup\mathcal{C}}G(I).\]
Finally, using subadditivity of $ V^* $ and the definition of
$\delta$, we have
\[V^*\left(\bigcup_{I\in\mathcal{B}\cup\mathcal{C}}G(I)\right)\le
\sum_{I\in\mathcal{B}\cup\mathcal{C}}V^*(G(I))=
\sum_{I\in\mathcal{B}}V^*(G(I))+\sum_{I\in\mathcal{C}}V^*(G(I))\le\]
\[\sum_{I\in\mathcal{B}}\frac{\ee}{2V(Q)}V(I)+L^m\sum_{I\in\mathcal{C}}V(I)
\le\frac{\ee}{2}+\frac{\ee}{2}=\ee.\]
}

\begin{lemma}
\label{fundamthm} 1. Let $ Q $ be a cube, denote by $ \mathcal{C}
$ the set of subcubes of $ Q $ and let $ \Phi:\mathcal{C}\to\R$ be
an additive Lipschitz function such that $ \Phi'(u)=0 $ holds for
almost all interior points $u$ of $ Q $. Then  $ \Phi $ is the
constant $ 0 $ function. 2. If $X\in\J$ and $\Psi\colon\J_X\to\R$
is an additive Lipschitz function such that $\Psi'(u)=0$ holds for
almost all interior points of $X$, then $\Psi$ is the constant $0$
function.
\end{lemma}

\proof{ 1. Suppose the contrary, then there exists a subcube $K$ such
that $\varepsilon:=|\Phi(K)|>0$. From the assumptions we have a
positive $L$ such that for each $I\in\mathcal{C}$, $|\Phi(I)|\le
L\cdot V(I)$, and we have a subset $H$ of $K$ with Lebesgue
measure $0$, which contains all the boundary points of $K$, such
that for all points $u\in K\setminus H$, $\Phi'(u)=0$.
Consequently, we have a countable set $\mathcal{T}$ of open
intervals with $\sum_{J\in\mathcal{T}}V(J)<\varepsilon/2L$ the
union of which covers $H$. To apply Cousin's lemma, define a
positive valued function $\delta$ on $K$. Assign to each $u\in H$
a $J_u\in \mathcal{T}$ that contains the point $u$ and then a
positive $\delta(u)$ such that $B(u,\delta(u))\subset J_u$, while
to each $u\in K\setminus H$, a $\delta(u)>0$ for which the
following implication holds: if a cube $I\in\mathcal{C}$ satisfies
the condition $u\in I\subset B(u,\delta(u))$ then
\[\left|\frac{\Phi(I)}{V(I)}\right|<\frac{\varepsilon}{2V(K)}.\]
Fix a $\delta$-fine dotted partition $\eta$ of $K$. The domain of
$\eta$ can be written as $\mathcal{A}\cup\mathcal{B}$ where for
$I\in \mathcal{A}$ and for $I\in \mathcal{B}$ we have $\eta_I\in
H$ and $\eta_I\in K\setminus H$, respectively. We get a
contradiction in the form $\varepsilon<\varepsilon$:
\[\varepsilon=|\Phi(K)|\le\sum_{I\in\mathcal{A}}|\Phi(I)|+
\sum_{I\in\mathcal{B}}|\Phi(I)|\le\]
\begin{equation}
\sum_{I\in\mathcal{A}}L\cdot V(I)+
\sum_{I\in\mathcal{B}}\frac{\ee\cdot V(I)}{2V(K)}\le
L\cdot\sum_{I\in\mathcal{A}}V(I)+\frac{\ee}{2}\stackrel{(*)}{<}
\frac{\ee}{2}+\frac{\ee}{2}=\ee.
\end{equation}
Inequality $(*)$ can be proved as follows. Using $\delta$-fineness
of $\eta$, each $I\in\mathcal{A}$ is a subset of $J_{\eta_I}$, so
the sum of volumes $V(I)$ for cubes $I$ belonging to the same
$J_u$ can be majorized by the volume of this common $J_u$,
consequently, for some finite subset $\mathcal{T}_0$ of
$\mathcal{T}$ we have
\[\sum_{I\in\mathcal{A}}V(I)\le\sum_{J\in\mathcal{T}_0}V(J)\le
\sum_{J\in\mathcal{T}}V(J) <\frac{\varepsilon}{2L}.\]
2. If
$Y\subset X$ is Jordan measurable, $\ee$ is a positive number and
$L>0$ is a Lipschitz constant for $\Psi$, then there exists a set
$H\subset\ir Y$ which is a finite union of cubes with
$V(Y\setminus H)<\ee/L$, therefore part 1.~of the theorem yields
\[|\Psi(Y)|=|\Psi(Y\setminus H)+\Psi(H)|=|\Psi(Y\setminus
H)|\le\ee.\]
}
\begin{remark}
Repeating a part of the proof of assertion 1.~we can get an
elementary proof of the fact that a bounded function $f$ defined
on a cube $K$ which is continuous in almost all interior points of
$K$, is integrable (consequently the same is true for a bounded
function defined on a Jordan measurable set (this is essentially
the same proof that one can find in \cite[Theorem 4]{Botsko:89}
for the case
$m=1$) : the domain is again a union $A\cup B$ where $A$ is of
small content and $B$ is a finite union of cubes). Indeed, let $
\varepsilon $ be a positive number. For any subcube $ I\subset K$,
define $L:=\osc_f(K)$, $H:=\partial K\cup\mbox{dis\,}f$. For
$ u\in H $, let the
definition of $\delta(u)$ be the same as in the previous proof,
while for $ u\in K\setminus H $, let $\delta(u)$ be any positive
number satisfying the condition $
\osc_f(\ol{B}(u,\delta(u)))<\ee/2V(K) $, and let the definitions
of $ \eta $, $\mathcal{A}$ and $\mathcal{B}$ be the same as
before. Then a proof of the inequality
\[\sum_{I\in\mathcal{A}}\osc_f(I)V(I)\,+\,\sum_{I\in\mathcal{B}}\osc_f(I)V(I)<
\varepsilon,\]
which implies integrability of $f$, can be formulated as follows.
Write the left hand side of this inequality followed by a ``$ \le
$'' sign and then switch to  line (1) and copy the previous
proof.\end{remark}

As three important corollaries, we give a characterization of the
`indefinite integral' of a given integrable function, a
characterization of the density functions of the constant zero
set-functions and a characterization of the set-functions
$\Psi\colon\J_X\to\R$ possessing a density function.

\begin{theorem}\label{charofintg}
If $X$ is a Jordan measurable set, $g:X\to\R$ an integrable
function and $ \Psi\colon\J_X\to\R$ then the following two
statements are equivalent: 1. $\Psi$ is an additive Lipschitz
function such that $ \Psi'(u)=g(u) $ holds for almost all interior
points $u$ of $X$, 2. $ \Psi(H)=\int_Hg $ for each $ H\in \J_X$
(in other words: $g$ is a density function of $\Psi$).
\end{theorem}

\proof{ 1.\bol2. Assertion 2.~of Lemma \ref{fundamthm} can be applied
to the set-function $ \J_X\ni H\mapsto \Psi(H)-\int_Hg$.

2.\bol1. Additivity and Lipschitz condition are well-known,
$\Psi'(u)=g(u)$ holds in each continuity points $u\in\ir X$ of
$g$.
}

\begin{theorem}\label{densofzero}
If $X$ is a Jordan measurable set, $g:X\to\R$ an integrable
function and $ \Psi\colon\J_X\to\R$ is the constant zero
set-function then the following three statements are equivalent:
1. $g(u)=0 $ holds for almost all points $u\in\ir X$, 2. $g$ is a
density function of $\Psi$, 3. $g(u)=0 $ holds for all continuity
points $u\in\ir X$ of $g$.
\end{theorem}

\proof{
1.\bol2. See assertion 1.\bol2. of Theorem \ref{charofintg}.
2.\bol3. See Remark \ref{indef}. 3.\bol1. Integrability implies
continuity in almost all interior points.
}

\begin{theorem}\label{charofindef}

Given a Jordan measurable set $X\subset\R^m$ and a set-function
$\Psi\colon\J_X\to\R$, the following two assertions are
equivalent: 1. $\Psi$ is an additive Lipschitz function that is
strongly differentiable in almost all interior points of $X$, 2.
$\Psi$ has a density function.
\end{theorem}
\proof{ 1.\bol2. We prove that the function $g\colon X\to \R$ defined
by
\[g(x):=\left\{
   \begin{array}{rl}
 \inf\{\sup\{\frac{\Psi(I)}{V(I)}\,:\,I\mbox{ is a subcube of }B(x,r)\}\;:\,r>0\},&\mbox{ if\ }x\in\ir  X,\\
0,&\mbox{ if\ }x\in \partial X\cap X
\end{array}\right.
\]
is a density function of $\Psi$. First we show that $g$ is
integrable, that is bounded, and continuous in almost all points
of $\ir X$. If $L$ is a Lipschitz constant for $\Psi$, then the
range of $g$ is contained in the interval $[-L,L]$, thus it
suffices to show that if $\Psi$ is strongly differentiable at an
interior point $u$ of $X$, then $g$ is continuous at $u$. Let
$u\in\ir X$, from the definition of strong differentiability we
have $\Psi'(u)=g(u)$. Let $\ee$ be a positive number and $\dd>0$
such that $B(u,\dd)\subset\ir X$ and for each subcube $I$ of
$B(u,\dd)$ $|\Psi(I)/V(I)-g(u)|<\ee$ holds. This implies that if
$\|x-u\|<\dd$ and $r<\dd-\|x-u\|$, then -- beeing each subcube $I$
of $B(x,r)$ a subcube of $B(u,\dd)$ --
\[\sup\left\{\frac{\Psi(I)}{V(I)}\,:\,I\mbox{ is a subcube of }B(x,r)\right\}\in [g(u)-\ee,g(u)+\ee],\]
therefore $g(x)\in[g(u)-\ee,g(u)+\ee]$, whenever $x\in B(u,\dd)$.
Now, Theorem\,\ref{charofintg} implies that $g$ is a density
function of $\Psi$.

2.\bol1. It is well-known that $\Psi$ is additive, Lipschitz, and
strongly differentiable in the continuity points $u\in\ir X$ of
$g$.
}

\section{Back to the change of variables}

\begin{theorem} \label{imageJoLi}
If $ A\subset\R^m $ is Jordan measurable, $ G:A\to\R^m $ is a
Lipschitz map and $ G $ is strongly differentiable at almost all
interior points of $ A $, then $ G(A) $ is Jordan measurable.
\end{theorem}

\proof{ Let $ L $ be a Lipschitz constant for $ G $, $ T $ the set of
those interior points of $ A $ where $ G $ is not strongly
differentiable, $ R $ and $ S $ the set of those points $x\in(\ir
A)\setminus T $, for which $ J_G(x) $ is regular or singular,
respectively. Finally, let $F$ be the unique continuous extension
of $G$ defined on $ \ol{A}$ (which is again a Lipschitz function
with Lipschitz constant $ L $). $ G(A) $ is bounded, because it is
a subset of the compact set $ F(\ol{A}) $. Continuity of $ F $
implies $ \ol{G(A)}=\ol{F(A)}=F(\ol{A}) $, therefore
\[\partial G(A)=\ol{G(A)}\setminus\ir G(A)=
F(\ol{A})\setminus\ir G(A)\]
\[=\left[G(R)\cup G(S)\cup G(T)\cup F(\partial A)\right]\setminus\ir G(A)\]
\[\subset\left[G(R)\setminus\ir G(A)\right]\cup G(S)\cup G(T)\cup F(\partial A)=
G(S)\cup G(T)\cup F(\partial A).\] The last equality follows from
the inclusion $ G(R)\subset\ir G(A)$ which is a consequence of
Theorem \ref{invfunthe}. Theorem \ref{genofSard} can be applied to
the function $ G\vert_{\ir A} $, from this we get that $ G(S) $ is
a Lebesgue-$0$-set. As both $ T $ and $
\partial A $ are Lebesgue-$0$-sets, to finish the proof it is
enough to observe that the image under a Lipschitz map of a
Lebesgue-$ 0 $-set is a Lebesgue-$ 0 $-set.
}

\begin{theorem}\label{nonoverl}
If $X\in\J$, $K\subset X$ is a set of Lebesgue measure $0$ and
$G\colon X\to \R$ is a Lipschitz function which is injective on
$\ir X\setminus K$, then for any two non-overlapping $A\in\J_X$
$B\in\J_X$, their images under $G$ are also non-overlapping.
\end{theorem}

\proof{ The inclusion
\[G(A)\cap G(B)\subset G(K)\cup G(\partial A\cap A)\cup
G(\partial B\cap B)\] follows from the fact that if $y=G(a)=G(b)$,
$a\in A\setminus K$ and $b\in B\setminus K$, then the relations
$a\in\ir A$, $b\in\ir B$ cannot hold at the same time: in the case
$a=b$ this would contradict to the fact that $A$ and $B$ are
non-overlapping, in the case $a\neq b$ -- to the injectivity
assumption. This inclusion implies that $G(A)\cap G(B)$ is of
Lebesgue measure $0$, therefore
\[[\ir G(A)]\cap[\ir G(B)]=\ir [G(A)\cap G(B)]=\emptyset.\]
}

All the existing proofs of the change of variables formula use the
following lemma, that we will also do.

\begin{lemma}
If the affine map $\ \ell:\R^m\to\R^m\,$ is defined by $\
\ell(x):=Ax+b$ where $A\in\R^{m\times m}$ and $b\in\R^m$, then for
each cube (in fact for each Jordan measurable set) $Q$ we have
$V(\ell(Q))=|\det A|\cdot V(Q)$.
\end{lemma}

\begin{theorem}\label{dibility}
Let $Q$ be a cube in $\R^m$, $G:Q\to\R^m$ a Lipschitz function
which is strongly differentiable at
$u\in\ir Q$ and $\Phi$ the cube-function defined on the set of
subcubes of $Q$ by $\Phi(I):=V(G(I))$.
Then $\Phi$ is strongly differentiable at $u$ and $\Phi'(u)=|\det J_G(u)|$.
\end{theorem}

\proof{
Set $D:=|\det J_G(u)|$. Suppose that $ D\neq0$ as the
other case has already been settled in Lemma \ref{lemmaSard}.
Let $\omega$ be a positive number and $\varepsilon\in(0,1)$ such that
\[D-\omega<(1-\ee)^mD\le(1+\ee)^mD<D+\omega.\]
Theorem \ref{invfunthe} yields a $\delta$ for this $\ee$; we may
and do suppose that $\ol{B}(u,\delta)\subset Q$. For each subcube
$I=\ol{B}(x,r)$ of $\ol{B}(u,\delta)$ we use the notations
\[I_-:=\ol{B}(x,(1-\ee)r),\qquad I_+:=\ol{B}(x,(1+\ee)r)\]
and apply Theorem \ref{invfunthe}:
\[D-\omega<(1-\ee)^mD=(1-\ee)^m\frac{V(\ell_x(I))}{V(I)}=
\frac{V(\ell_x(I_-))}{V(I)}\le\frac{V(G(I))}{V(I)}\le\]
\[\frac{V(\ell_x(I_+))}{V(I)}\le(1+\ee)^m\frac{V(\ell_x(I))}{V(I)}\le
(1+\ee)^mD<D+\omega.\]
}

\begin{theorem}\label{konkrsfv}
If $X\subset\R^m$ is a Jordan measurable set, $K\subset X$ a
Lebesgue-$0$-set, $G:X\to\R^m$ a Lipschitz map which is strongly
differentiable in almost all points of $\,\ir X$, injective on
$X\setminus K$ and the set-function $\J_X\ni H\mapsto V(G(H))$ is
denoted by $\Psi$, then a density function $g$ of $\Psi$ can be
constructed in this way: using the notation
\[B^j(x,r):=\{(y,z)\in B(x,r)\times B(x,r):\,y_j\neq z_j, \
y_i=z_i\ \mbox{if}\ \ i\neq j\}\ \  (x\in\R^m, \ r>0),\] let
$g\colon X\to\R$ be the function $x\mapsto|\det(\ol{g}_{ij}(x))|$,
where for all pairs of integers $1\le i,j\le m$, the definition of
the functions $\ol{g}_{ij}\colon X\to \R$ is
\[\ol{g}_{ij}(x):=\left\{
   \begin{array}{rl}
 \inf\{\sup\{\frac{(G_i(z)-G_i(y)}{z_j-y_j}\,:\,(y,z)\in B^j(x,r)\}\;:\,r>0\},&\ \mbox{ if }\ x\in\ir X,\\
0,&\ \mbox{ if }\ x\in \partial X\cap X
\end{array}\right..\]
\end{theorem}
\proof{ First we show that for each pair $(i,j)$, $\ol{g}_{ij}$ is
integrable, that is bounded and almost everywhere in $\ir X$
continuous. If $L$ is a Lipschitz constant for $G$, then the range
of $\ol{g}_{ij}$ is contained in the interval $[-L,L]$, thus it
suffices to show that if $G$ is strongly differentiable at an
interior point $u$ of $X$, then $\ol{g}_{ij}$ is continuous at
$u$. Let $u\in\ir X$ such that $G$ is strongly differentiable at
$u$ and let $\ee$ be a positive number. From Remarks \ref{edked}
and \ref{edepd} we know that $G_i$ is strongly partially
differentiable at $u$ with respect to the $j$-th variable and from
Definition\,\ref{epd} it is clear that
$\ol{g}_{ij}(u)=D^s_jG_i(u)$. Select a $\dd>0$ such that
\[\frac{(G_i(z)-G_i(y)}{z_j-y_j}\in[\ol{g}_{ij}(u)-\ee,\ol{g}_{ij}(u)+\ee],\quad\mbox{whenever}\quad
(y,z)\in B^j(x,\dd),\]
If $\|x-u\|<\dd$ and $r<\dd-\|x-u\|$, then
-- being $B^j(x,r)\subset B^j(u,\dd)$ --
\[\sup\left\{\frac{G_i(z)-G_i(y)}{z_j-y_j}\,:\,(y,z)\in B^j(x,r)\right\}\in [\ol{g_{ij}}(u)-\ee,\ol{g}_{ij}(u)+\ee],\]
therefore
$\ol{g}_{ij}(x)\in[\ol{g_{ij}}(u)-\ee,\ol{g}_{ij}(u)+\ee]$,
whenever $x\in B(u,\dd)$. Now, using continuity of the function
$\det\colon\R^{m\times m}\to\R$ and compactness of $[-L,L]^{m^2}$
in $\R^{m^2}$ we get that $g$ is bounded, and continuous at the
strong differentiability points of $G$, thus $g$ is integrable.
Moreover, at the strong differentiability points $u\in\ir X$ of
$G$, we have $\Psi'(u)=g(u)$ (see Theorem\,\ref{dibility}) and
according to Theorems \ref{nonoverl}, \ref{VGHLip}, $\Psi$ is an
additive Lipschitz function. These facts and
Theorem\,\ref{charofintg} imply that $g$ is a density function of
$\Psi$.
}

\begin{theorem}\label{maintheor}
If $X\subset\R^m$ is a Jordan measurable set, $K\subset X$ a
Lebesgue-$0$-set, $G:X\to\R^m$ a Lipschitz map which is strongly
differentiable in almost all points of $\ir X$, injective on
$X\setminus K$, and $f:G(X)\to\R$ is any bounded function, then

\noindent 1. there exists an integrable function
$\ol{g}:X\to\R^{m\times m}$ such that for almost all $x\in\ir X$,
$\ol{g}(x)=J_G(x)$,

\noindent 2. for each integrable function $\ol{h}:X\to\R^{m\times
m}$ with this property, the function
\[\psi:X\to\R,\quad \psi(x):=f(G(x))\cdot|\det \ol{h}(x)|\]
is integrable if and only if $f$ is integrable, and
\begin{equation} \int_{G(X)}f=\int_X\psi\end{equation}
holds whenever one of $f$ and $\psi$ is (that is both of $f$ and
$\psi$ are) integrable.
\end{theorem}

\proof{
1. From Theorem \ref{konkrsfv} we already know that the function
$\ol{g}$ defined there is integrable, and in the strong
differentiability points $x$ of $G$, $\ol{g}(x)=J_G(x)$ holds.

\noindent 2. Theorem \ref{alap} can be applied with $X\ni
x\mapsto|\det\ol{h}(x)|=:g(x)$. Indeed, Theorem \ref{imageJoLi}
implies that condition a) of Theorem \ref{alap} is satisfied,
while b) is implied by Theorem \ref{nonoverl}. As for condition
c), $g$ differs from the function $X\ni x\mapsto|\det\ol{g}(x)|$
on a set of Lebesgue measure $0$; according to Theorem
\ref{konkrsfv}, the latter is a density function of $\J_X\ni
H\mapsto V(G(H))$, so the same is true for $g$ (see Theorem
\ref{densofzero}).
}
\begin{remark}
In Theorem \ref{maintheor}, the injectivity assumption cannot be
omitted, in particular, Theorem~24.26 in \cite{Bartle:64} is
false. Counterexample: $m=2$, $A:=[1,2]\times[0,4\pi]$,
$G(x,y):=(e^x\cos y,e^x\sin y)$, $f(x,y):=1$.
\end{remark}

\bibliographystyle{plain}

\end{document}